\documentclass[11pt,oneside,reqno]{amsart}

\usepackage{amssymb}
\usepackage{algpseudocode}
\usepackage{algorithm}
\usepackage{subcaption}
\usepackage{tikz}
\usepackage[singlelinecheck=false]{caption}
\usepackage{verbatim}
\usepackage{url}
\usepackage{graphicx}
\usepackage[noadjust]{cite}
\RequirePackage{dsfont} 
 \scrollmode
\allowdisplaybreaks
\usepackage{graphicx}
\usepackage{epstopdf}
\usepackage{bbm}
\usepackage[utf8]{inputenc}
\usepackage[T1]{fontenc}

\usepackage{hyperref}

\makeatletter
\newcommand{\descitem}[2]{
  \item[#1]
  \phantomsection
  \def\@currentlabel{#1}
  \label{#2}
}
\makeatother

\definecolor{red-new}{rgb}{1,0,0}
\definecolor{dark-blue}{rgb}{0,0,1}
\definecolor{light-blue}{rgb}{0,0.6,0.6}
\definecolor{gray-new}{rgb}{0.3,0.3,0.3}

\usepackage{amsmath}
\usepackage{tikz}

\makeatletter
\newcommand{\xleftrightarrow}[2][]{\ext@arrow 3359\leftrightarrowfill@{#1}{#2}}
\makeatother

\newcommand{\xdasharrow}[2][->]{
\tikz[baseline=-\the\dimexpr\fontdimen22\textfont2\relax]{
\node[anchor=south,font=\scriptsize, inner ysep=1.5pt,outer xsep=2.2pt](x){#2};
\draw[shorten <=3.4pt,shorten >=3.4pt,dashed,#1](x.south west)--(x.south east);
}
}

\newcommand{\DEBUG}{}

\ifdefined\DEBUG

  \def\rem#1{{\marginpar{\raggedright\scriptsize #1}}}
  \newcommand{\pmr}[1]{\rem{\color{blue}{$\bullet$ #1}}}
  \newcommand{\ppr}[1]{\rem{\color{red}{$\bullet$ #1}}}
 \else

  \newcommand{\ppr}[1]{}
  \newcommand{\pmr}[1]{}
 \fi

\theoremstyle{plain}
\newtheorem{theorem}{Theorem}
\newtheorem{lemma}{Lemma}
\newtheorem{fact}{Fact}
\newtheorem{corollary}{Corollary}

\theoremstyle{definition}
\newtheorem{remark}{Remark}

\usepackage{enumitem}
\begin{document}

\title
[Randomized Euler scheme for nonlinear ODEs]
{On the performance of the randomized Euler scheme for nonlinear ODEs under inexact information}

\author[T. Bochacik]{Tomasz Bochacik}
\address{AGH University of Krakow,
Faculty of Applied Mathematics,\newline
Al. A.~Mickiewicza 30, 30-059 Krak\'ow, Poland}
\email{bochacik@agh.edu.pl, corresponding author}

\author[B. Fober]{Bartosz Fober}
\address{HSBC Service Delivery (Polska) Sp. z o.o.,\newline
ul. Kapelanka 42A, 30-347 Krak\'ow, Poland}
\email{bartosz.fober@hsbc.com}

\author[P. Przyby{\l}owicz]{Pawe{\l} Przyby{\l}owicz}
\address{AGH University of Krakow,
Faculty of Applied Mathematics,\newline
Al. A.~Mickiewicza 30, 30-059 Krak\'ow, Poland}
\email{pprzybyl@agh.edu.pl}

\begin{abstract}
We analyze the randomized explicit Euler scheme for nonlinear ordinary differential equations under a local one-sided Lipschitz condition, local H\"older regularity in the space variable, and H\"older regularity in the time variable. We allow for noisy evaluations of the right-hand side. For the resulting approximation, we derive an error bound in which the contribution associated with the H\"older regularity in the time variable exhibits the characteristic $1/2$-order gain of randomized methods. While existing results cover deterministic schemes under the one-sided Lipschitz condition and randomized schemes under the classical Lipschitz condition, our results combine randomized approximation and one-sided Lipschitz regularity.
\end{abstract}
\maketitle
\tableofcontents

Numerical approximation of solutions of ordinary differential equations (ODEs) is a fundamental task in scientific computing, with applications ranging from mechanics and control theory to biological, chemical, and engineering models. In this paper, we consider the randomized explicit Euler scheme for a certain class of nonlinear ODEs. In particular, we assume that the right-hand side satisfies the local one-sided Lipschitz condition and the local H\"older condition with respect to the space variable. For the time variable, we assume H\"older continuity. We additionally consider the practically important case in which every evaluation of the right-hand side is contaminated by additive noise satisfying an at-most-linear growth condition with respect to the space variable. Our main result provides an $L^p(\Omega)$-error bound for the randomized explicit Euler scheme under these assumptions. The bound reveals the characteristic benefit of randomization: the term resulting from the H\"older regularity in the time variable contains an additional factor corresponding to a gain of $1/2$ in the convergence rate.

This research is primarily motivated by the analysis of the deterministic Euler scheme in \cite{CzyzMorPrzyb}. In that work, the same classes of right-hand side and noise functions were considered. The error bound obtained therein provides a natural benchmark for the randomized scheme studied in the present paper.

The assumptions considered in \cite{CzyzMorPrzyb} and in the present paper are motivated by two features that are relevant in practical computations. First, the classical global Lipschitz condition can be restrictive for nonlinear problems. In our setting, we weaken this assumption in two respects: we require only a one-sided Lipschitz condition, and we impose it locally, so that the corresponding Lipschitz constant is allowed to depend on the compact subset rather than being uniform over the whole space. Such an assumption covers a broader class of nonlinear problems while still allowing error bounds for their numerical approximation. Second, in practical applications, the right-hand side is generally not available exactly. Its evaluations may be affected by measurement errors, finite-precision effects, errors in the computation of model quantities, or uncertainty in model parameters. This motivates the analysis of numerical methods based on noisy information, where evaluations of the right-hand side are perturbed. Error bounds in this setting are also important because they allow one to separate the asymptotic benefit of increasing the number of time steps from the error floor imposed by the accuracy of the available information.

Randomization has emerged as an effective approach to the numerical solution of both deterministic ODEs and stochastic differential equations (see, e.g., \cite{PWX, BKMPS}), particularly for problems with limited regularity. The underlying idea is to replace deterministic evaluation points with randomly chosen points, which can lead to the cancellation of local errors and, consequently, to improved convergence rates. 

For ODEs, a randomized approach has been developed in a number of works, both from the perspective of information-based complexity \cite{Kac1, HeinMilla, daun1} and for specific randomized schemes \cite{KruseWu_1, backward_euler, JenNeuen, BPS}. In particular, these results demonstrate that randomization can improve the convergence rate by $1/2$ compared with deterministic schemes under suitable regularity assumptions. Furthermore, randomized Euler and two-stage Runge--Kutta schemes have been analyzed under inexact information in \cite{randEuler, randRK}. The present work builds on these developments by considering a randomized explicit Euler scheme and incorporating noisy evaluations of the right-hand side. The main contribution of this paper is that the error analysis is carried out under substantially weaker regularity assumptions, which require different proof techniques from those used in the existing analyses.

\section{Preliminaries} \label{sec:prel}

\subsection{Initial value problem}

Consider initial value problems of the following form: 
\begin{equation}
		\left\{ \begin{array}{ll}
			z'(t)= f(t,z(t)), \ t\in [a,b], \\[2pt]
			z(a) = \eta, 
		\end{array}\right.
        \label{eq:ivp}
\end{equation}
where $-\infty < a < b < \infty$, $\eta\in\mathbb{R}^d$, $f\colon [a,b]\times\mathbb{R}^d\to\mathbb{R}^d$, $d\in\mathbb{Z}_+$. 

\subsection{Randomized scheme}

Let $n\in\mathbb{Z}_+$ and $h=\frac{b-a}{n}$, $t_j = a+jh$ for $j\in\{0,1,\ldots,n\}$. Moreover, let $\tau_1,\tau_2,\ldots$ be independent random variables uniformly distributed on the interval $[0,1]$, defined on a common probability space $(\Omega,\Sigma,\mathbb{P})$, where the $\sigma$-algebra $\Sigma$ is complete. We also consider a filtration $(\mathcal{F}_j)_{j=0}^\infty$ generated by random variables $\tau_1,\tau_2,\ldots$. Specifically, $\mathcal{F}_0 = \sigma(N)$, where $N=\{ A\in\Sigma\colon \mathbb{P}(A)=0\}$, and $\mathcal{F}_j = \sigma\bigl( \sigma(\tau_1,\ldots,\tau_j), N \bigr)$. For $j\in\{1,\ldots,n\}$, let $\theta_j = t_{j-1} + h\tau_j$. We use the notation $\langle\cdot,\cdot\rangle$ to denote the standard scalar product in $\mathbb{R}^d$, while the Euclidean norm on $\mathbb{R}^d$, induced by the above scalar product, is denoted by $\|\cdot\|$.

We investigate the randomized explicit Euler scheme. The approximate value of $z(t_j)$, where $z$ is the exact solution of \eqref{eq:ivp}, will be denoted by $Y_j$. The sequence $(Y_j)_{j=0}^n$ is obtained recursively:
\begin{equation}
    \left\{ \begin{array}{ll}
         Y_0 = \eta, &  \\
         Y_j = Y_{j-1} + h\cdot f(\theta_j, Y_{j-1}) & \text{for } j\in\{1,\ldots,n\}.
    \end{array} \right.
    \label{eq:rand_ee}
\end{equation}
By induction, $Y_j\in \mathcal{F}_j$ for $j\in\{0,1,\ldots,n\}$.

Error bounds for this scheme were established in \cite{KruseWu_1, randEuler} under the classical Lipschitz condition. In this paper, we aim to further weaken assumptions imposed on the right-hand side function $f$, in particular by considering the local one-sided Lipschitz condition.

\section{Error bounds} \label{sec:error}

\subsection{Class of initial value problems}

We assume that the right-hand side function $f$ in the equation \eqref{eq:ivp} satisfies the following conditions:
\begin{description}
\setlength{\itemindent}{-0.28in}
\setlength\itemsep{4pt}
\descitem{G1}{itm:G1} $f\in \mathcal{C}\left([a,b]\times\mathbb{R}^d, \mathbb{R}^d\right)$.
\descitem{G2}{itm:G2} There exists $K>0$ such that $$\| f(t,x)\| \leq K\left(1+\left\|x\right\|\right)$$ for all $(t,x)\in [a,b]\times\mathbb{R}^d$.
\descitem{G3}{itm:G3} For each compact set $U\subset \mathbb{R}^d$ there is a constant $H_U\in\mathbb{R}$ such that $$\langle x_1-x_2, f(t,x_1)-f(t,x_2)\rangle \leq H_U \|x_1-x_2\|^2$$ for all $t\in[a,b]$ and $x_1,x_2\in U$.
\descitem{G4}{itm:G4} There exist $\alpha,\beta\in (0,1]$ such that for every compact set $U\subset \mathbb{R}^d$ there is a constant $L_U>0$ such that $$\|f(t_1,x_1)-f(t_2,x_2)\|\leq L_U\bigl( |t_1-t_2|^\alpha + \|x_1-x_2\|^\beta \bigr)$$
for all $t_1,t_2\in [a,b]$ and $x_1,x_2 \in U$.
\end{description}

\subsection{Noisy information}

Let $\delta\in[0,1]$ be a precision parameter. Consider an IVP \eqref{eq:ivp} with fixed initial condition $\eta$ and right-hand side function $f$. Let us assume that we have access only to a noisy value $\tilde\eta$ of $\eta$ such that:
\begin{equation}
    \|\tilde \eta -\eta\|\leq \delta
    \label{eq:noisy_eta}
\end{equation}
and to noisy evaluations of the function $f$:
\begin{equation}
    \tilde f(t,x) = f(t,x)+\delta_f(t,x) \ \text{ for } \ (t,x)\in [a,b]\times\mathbb{R}^d,
    \label{eq:noisy_f}
\end{equation}
where the noise function $\delta_f\colon [a,b]\times\mathbb{R}\to\mathbb{R}^d$ has at most linear growth in the space variable:
\begin{equation}
    \|\delta_f(t,x)\|\leq \delta (1+\|x\|)  \ \text{ for } \ (t,x)\in [a,b]\times\mathbb{R}^d.
    \label{eq:noise_function}
\end{equation}

Randomized Euler scheme \eqref{eq:rand_ee} based on noisy information takes the following form:
\begin{equation}
    \left\{ \begin{array}{ll}
         \tilde Y_0 = \tilde\eta, &  \\
         \tilde Y_j = \tilde Y_{j-1} + h\cdot\tilde f(\theta_j, \tilde Y_{j-1}) & \text{for } j\in\{1,\ldots,n\}.
    \end{array} \right.
    \label{eq:noisy_rand_ee}
\end{equation}

\subsection{Bounds on the exact, approximate, and local solutions}

The existence, uniqueness, and regularity of the solution to \eqref{eq:ivp} are ensured by the following lemma, proven in \cite{CzyzMorPrzyb} as Lemma A.2.

\begin{lemma} \label{lemma:exist_exact}
    Let $\eta\in\mathbb{R}^d$ and let us assume that function $f\colon [a,b]\times\mathbb{R}^d\to\mathbb{R}^d$ satisfies assumptions \eqref{itm:G1}--\eqref{itm:G3}. Then there exists a unique solution $z=z(\eta,f)\in\mathcal{C}^1([a,b], \mathbb{R}^d)$ of \eqref{eq:ivp}. Moreover, \begin{equation}
        \label{eq:sup_z_bound}
        \sup_{a\leq t \leq b} \|z(t)\| \leq K_0,
    \end{equation}
    where $K_0=K_0(a,b,K,\eta) = e^{K(b-a)}\bigl( \|\eta\| + K(b-a)\bigr)$.
\end{lemma}

The following two facts provide upper bounds of the norm of approximate solution and a local solution, which in a sense is a bridge between the exact and the approximate solutions.

\begin{fact} \label{lemma:bound_approx_sol}
    Let $\eta\in\mathbb{R}^d$, $f$ satisfy \eqref{itm:G2}, and $\delta\in [0,1]$. Then 
    $$\max_{n\in\mathbb{Z}_+} \max_{k\in\{0,1,\ldots,n\} }\|\tilde Y_k\|\leq e^{(K+1)(b-a)}(2+\|\eta\|) $$
    with probability $1$, where $\tilde Y_k$s (for a fixed $n$) are given by \eqref{eq:noisy_rand_ee}.
\end{fact}

\begin{proof}
    Fix $n\in\mathbb{Z}_+$. By \eqref{eq:rand_ee}, for every $k\in\{0,1\ldots,n-1\}$,
    \begin{align*}
        \|\tilde Y_{k+1}\| & \leq \|\tilde Y_{k}\| + h \|f(\theta_{k+1}, \tilde Y_k)\| +h \|\delta_f(\theta_{k+1}, \tilde Y_k)\| \leq \|\tilde Y_{k}\| + h(K+\delta)(1+\|\tilde Y_{k}\|) \\ & \leq (1+h(K+1))\|\tilde Y_k\|+h(K+1).
    \end{align*}
    The second inequality follows from \eqref{itm:G2} and \eqref{eq:noise_function}. By discrete Gronwall's inequality, for every $k\in\{0,1\ldots,n\}$,
    \begin{align*}
        \| \tilde Y_k\| & \leq (1+h(K+1))^k \|\tilde Y_0\| + (1+h(K+1))^k-1 \leq e^{(K+1)hn}\|\tilde \eta\|+e^{(K+1)hn}-1\\ & \leq e^{(K+1)(b-a)}(2+\|\eta\|)
    \end{align*}
   with probability $1$. This completes the proof.
\end{proof}

\begin{remark}
    If $\delta=0$, then $\tilde Y_k = Y_k$ for $k\in\{0,1,\ldots,n\}$. Thus, for the randomized Euler scheme based on exact information \eqref{eq:rand_ee}, we also have
    $$\max_{n\in\mathbb{Z}_+} \max_{k\in\{0,1,\ldots,n\} }\|Y_k\|\leq e^{(K+1)(b-a)}(2+\|\eta\|).$$
\end{remark}

\begin{fact} \label{lemma:local_ivp}
    For $n\in\mathbb{Z}_+$ and $k\in\{0,1,\ldots,n-1\}$ consider a local IVP:
    \begin{equation}
        \label{eq:local_ivp}
        \left\{ \begin{array}{ll}
			z_{k}'(t) = f(t, z_k(t)), \ t\in [t_k,t_{k+1}], \\[2pt]
			z_k(t_k) = Y_k.
		\end{array}\right.
    \end{equation}
    This IVP has a unique solution $z_k \colon [t_k,t_{k+1}]\times\Omega\to\mathbb{R}^d$. Moreover, $z_k\in \mathcal{C}^1 ([t_k,t_{k+1}], \mathbb{R}^d)$ and $z_k(t)$ for $t\in [t_k,t_{k+1}]$ is $\mathcal{F}_k$-measurable.  Finally, with probability $1$, for every $n\in\mathbb{Z}_+$, $k\in\{0,1,\ldots,n-1\}$, and $t\in[t_k,t_{k+1}]$,
    $$\|z_k(t)\| \leq e^{K(b-a)}\bigl( e^{K(b-a)} + K(b-a) \bigr)(1+\|\eta\|).$$
\end{fact}

\begin{proof} 
   The first part of the fact (related to the existence of the solution) follows from Lemma \ref{lemma:exist_exact}. For the second part (bound), let us fix $n\in\mathbb{Z}_+$, consider $k\in\{0,1,\ldots,n-1\}$ and $t\in [t_k,t_{k+1}]$, and observe that 
   \begin{align*}
       \|z_k(t)\| & = \Bigl\| Y_k + \int\limits_{t_k}^t f(s,z_k(s))\,\mathrm{d}s \Bigr\| \leq \|Y_k\| + \int\limits_{t_k}^t \|f(s,z_k(s))\|\,\mathrm{d}s \leq \|Y_k\| + \int\limits_{t_k}^t (1+\|z_k(s)\|)\,\mathrm{d}s \\ & \leq e^{K(b-a)}(1+\|\eta\|) + K(b-a) + K\int\limits_{t_k}^t \|z_k(s)\|\,\mathrm{d}s
   \end{align*}
   with probability $1$. In the last passage of the first line, we used \eqref{itm:G2}. In the passage to the second line, we utilized Fact \ref{lemma:bound_approx_sol}. By the continuous Gronwall's inequality,
   \begin{align*}
       \|z_k(t)\| & \leq e^{K(b-a)}\bigl( e^{K(b-a)}(1+\|\eta\|) + K(b-a) \bigr) \\ & \leq e^{K(b-a)}\bigl( e^{K(b-a)} + K(b-a) \bigr)(1+\|\eta\|),
   \end{align*}
   which completes the proof.
\end{proof}

\begin{corollary}
\label{cor:R}
By Lemma \ref{lemma:exist_exact}, Fact \ref{lemma:bound_approx_sol}, and Fact \ref{lemma:local_ivp}, with probability $1$,
$$z(t), z_k(t), Y_k, Y_n, \tilde Y_k, \tilde Y_n \in \bar B_0(R^*)$$
for every $n\in\mathbb{Z}_+$, $k\in\{0,1,\ldots,n-1\}$, and $t\in [a,b]$, where
\begin{equation}
    R^* = e^{(K+1)(b-a)}\bigl( e^{K(b-a)} + K(b-a) \bigr)(2+\|\eta\|)
    \label{eq:Rstar}
\end{equation}
and $\bar B_0(R^*) = \{x\in\mathbb{R}^d \colon \|x\|\leq R^*\}$.
\end{corollary}
\begin{proof}
 Consider upper bounds obtained for $\|z(t)\|$, $\|\tilde Y_k\|$, and $\|z_k(t)\|$ in Lemma \ref{lemma:exist_exact}, Fact \ref{lemma:bound_approx_sol}, and Fact \ref{lemma:local_ivp}, respectively.
\end{proof}

\subsection{Results on the $L^p(\Omega)$-error}

Consider the initial value problem \eqref{eq:ivp}, fix $n\in\mathbb{Z}_+$, and set $h=\frac{b-a}{n}$. Let $z\colon [a,b]\to\mathbb{R}^d$ be an exact solution to this IVP and let $(Y_k)_{k=0}^n$ be the sequence of approximated values of $z$ at grid points $(t_k)_{k=0}^n=(a+hk)_{k=0}^n$, generated by the randomized Euler scheme \eqref{eq:rand_ee}. Let $$e_k:= z(t_k)-Y_k$$ for $k\in\{0,1,\ldots,n\}$ and $$u_k:=z(t_k)-z_{k-1}(t_k)$$ for $k\in\{1,\ldots,n\}$. Recall that the local IVP is given by \eqref{eq:local_ivp}. 

By linear interpolation based on points $(t_k,\tilde Y_k)_{k=0}^n$ produced by \eqref{eq:noisy_rand_ee}, one obtains an approximate solution to \eqref{eq:ivp} defined on the entire interval $[a,b]$, which will be denoted as $\tilde y_n\colon [a,b]\to\mathbb{R}^d$. The explicit expression for $\tilde y_n$ is 
\begin{equation}
    \tilde y_n(t) = \tilde Y_{k-1} + \frac{t-t_{k-1}}{t_k-t_{k-1}}\cdot (\tilde Y_k-\tilde Y_{k-1}) \ \ \text{ for } \ \ t\in[t_{k-1},t_k], \ \ \ k\in\{1,\ldots,n\}.
    \label{eq:barz}
\end{equation}
Analogously, $y_n$ would be a piecewise linear function based on points $(t_k,Y_k)_{k=0}^n$ produced by the randomized Euler scheme under exact information \eqref{eq:rand_ee}.

Let us also consider a piecewise linear function based on the actual values of an exact solution $z$ to \eqref{eq:ivp}:
\begin{equation}
    \bar z_n (t) = z(t_{k-1}) + \frac{t-t_{k-1}}{t_k-t_{k-1}}\cdot \bigl(z(t_k)-z(t_{k-1})\bigr) \ \ \text{ for } \ \ t\in[t_{k-1},t_k], \ \ \ k\in\{1,\ldots,n\}.
    \label{eq:z-lin}
\end{equation}

\begin{fact}
    Let $n\in\mathbb{Z}_+$ and consider an IVP \eqref{eq:ivp} with function $f\colon [a,b]\times\mathbb{R}^d\to\mathbb{R}^d$ satisfying the conditions \eqref{itm:G1}--\eqref{itm:G4}. Then the error between the exact solution $z\colon [a,b]\to\mathbb{R}^d$ of the IVP and its piecewise linear approximation $\bar z_n$ given by \eqref{eq:z-lin} can be bounded as follows:
    \begin{equation}
        \sup_{t\in [a,b]} \|z(t)-\bar z_n(t)\| \leq L_1dh^{\alpha+1} +L_1dK^\beta(1+K_0)^\beta h^{\beta+1},
        \label{eq:bound-z-lin-z}
    \end{equation}
    where $L_1 = L_{\bar B_0(R^*)}$, cf. \eqref{itm:G4} and \eqref{eq:Rstar}, $K$ is a constant from \eqref{itm:G2}, and $K_0$ is defined in Lemma \ref{lemma:exist_exact}.
    \label{f:bound-z-lin-z}
\end{fact}

The fact below establishes an upper bound on the error of piecewise linear interpolation to the exact solution $z$ of the IVP \eqref{eq:ivp}.

\begin{proof}
    Fix $k\in\{1,\ldots,n\}$ and $t\in [t_{k-1},t_k]$. By Lagrange's mean value theorem, there exist constants $\alpha_1,\ldots,\alpha_d\in (t_{k-1},t_k)$ such that
    \begin{equation}
        z(t) = \bigl(z_1(t),\ldots,z_d(t)\bigr)= \bigl(z_1(t_{k-1})+z'_1(\alpha_1)(t-t_{k-1}),\ldots,z_d(t_{k-1})+z'_d(\alpha_d)(t-t_{k-1})\bigr)
        \label{eq:z-alpha}
    \end{equation}
    Moreover, for every $l\in\{1,\ldots,d\}$, there exists $\beta_l\in (t_{k-1},t_k)$ such that
    $$\frac{z_l(t_k)-z_l(t_{k-1})}{t_k-t_{k-1}} = z'_l(\beta_l).$$
    Inserting this into \eqref{eq:z-lin}, we obtain 
    $$\bar z_{n} (t) = \bigl( z_1(t_{k-1}) + z'_1(\beta_1)\cdot (t-t_{k-1}), \ldots, z_d(t_{k-1}) + z'_d(\beta_d)\cdot (t-t_{k-1}) \bigr).$$
    for $t\in [t_{k-1},t_k]$. Combining this with \eqref{eq:z-alpha},
    \begin{align*}
        \|z(t)-\bar z_n(t)\| & = \sum_{l=1}^d \bigl| z'_l(\alpha_l)-z'_l(\beta_l) \bigr|(t-t_{k-1}) \leq h \sum_{l=1}^d\|z'(\alpha_l)-z'(\beta_l)\| \\ & = h \sum_{l=1}^d\|f(\alpha_l,z(\alpha_l))-f(\beta_l,z(\beta_l))\| \leq L_1h \sum_{l=1}^d \bigl( |\alpha_l-\beta_l|^\alpha +\|z(\alpha_l)-z(\beta_l)\|^\beta \bigr) \\ & \leq  L_1h \sum_{l=1}^d\bigl(h^\alpha +\|z(\alpha_l)-z(\beta_l)\|^\beta\bigr),
    \end{align*}
    where in the second line, we used \eqref{itm:G4} with $U=\bar B_0(R^*)$, cf. Corollary \ref{cor:R}, that is,
    \begin{equation}
        L_1 = L_{\bar B_0(R^*)}.
        \label{eq:L1}
    \end{equation}
   Moreover, by \eqref{itm:G2} and \eqref{eq:sup_z_bound} from Lemma \ref{lemma:exist_exact}, 
    \begin{equation*}
        \|z(t)-z(s)\| \leq \int\limits_s^t \|z'(u)\|\,\mathrm{d}u = \int\limits_s^t \|f(u,z(u))\|\,\mathrm{d}u \leq K(1+K_0)(t-s)
    \end{equation*}
    for $a\leq s<t\leq b$. Hence,
    $$\|z(t)-\bar z_n(t)\| \leq L_1h \sum_{l=1}^d\bigl(h^\alpha +K^\beta(1+K_0)^\beta|\alpha_l-\beta_l|^\beta\bigr)\leq L_1dh^{\alpha+1} +L_1dK^\beta(1+K_0)^\beta h^{\beta+1},$$
    which completes the proof.
\end{proof}

An error due to informational noise is bounded in the lemma below. Proof technique was originally established in \cite{CzyzMorPrzyb} (see Theorem 4.1).

\begin{lemma}
    \label{lemma:Y-tilde}
    For fixed $\eta\in\mathbb{R}^d$, $f\colon [a,b]\times \mathbb{R}^d\to\mathbb{R}^d$ satisfying \eqref{itm:G1}--\eqref{itm:G4}, let $(Y_k)_{k=0}^n$ and $(\tilde Y_k)_{k=0}^n$ be approximate solutions to \eqref{eq:ivp} generated by schemes \eqref{eq:rand_ee} and \eqref{eq:noisy_rand_ee}, respectively. Then, there exists a constant $\tilde K>0$ such that for every $n\in\mathbb{Z}_+$ satisfying $n\geq 2(b-a)$,
    \begin{equation}
        \max_{k\in\{0,1,\ldots,n\}}\|\tilde Y_k-Y_k\|\leq \tilde K \bigl( \delta + h^{\frac12}\cdot\mathbbm{1}_{(0,1]}(\delta)\cdot\mathbbm{1}_{(0,1)}(\beta) \bigr)
        \label{eq:max-Y-tilde-Y}
    \end{equation}
    with probability $1$.
\end{lemma}
\begin{proof}
    Let $\Delta_k = \tilde Y_k-Y_k$ for $k\in \{0,1,\ldots,n\}$. For $k=0$, by \eqref{eq:noisy_eta},
    $$\|\Delta_0\|=\|\tilde \eta -\eta\|\leq \delta.$$
    Now, let $k\in\{1,\ldots,n\}$. Then
    \begin{equation}
        \Delta_k = \Delta_{k-1} + h\bigl(f(\theta_k,\tilde Y_{k-1})-f(\theta_k,Y_{k-1})\bigr)+h\delta_f(\theta_k,\tilde Y_{k-1})
        \label{eq:delta-k}
    \end{equation}
    For $\beta=1$, by \eqref{itm:G4}, \eqref{eq:noise_function}, \eqref{eq:L1}, and then by Corollary \ref{cor:R},
    \begin{align*}
        \|\Delta_k\| & \leq \|\Delta_{k-1}\| + h\bigl\|f(\theta_k,\tilde Y_{k-1})-f(\theta_k,Y_{k-1})\bigr\|+h\|\delta_f(\theta_k,\tilde Y_{k-1})\| \\ & \leq (1+ L_1h)\|\Delta_{k-1}\| + \delta h (1+\|\tilde Y_{k-1}\|) \leq (1+ L_1h)\|\Delta_{k-1}\| + \delta h (1+R^*)
    \end{align*}
    for $k\in\{1,\ldots,n\}$. Thus, by \eqref{eq:noisy_eta},
    \begin{align*}
        \|\Delta_k\| & \leq (1+ L_1h)^k\|\Delta_{0}\| + \delta h (1+R^*) \cdot \frac{(1+L_1h)^k-1}{L_1h} \leq e^{L_1(b-a)}\Bigl(1+\frac{1+R^*}{L_1}\Bigr)\delta.
    \end{align*}
    Now we consider the case $\beta\in (0,1)$. Let $k\in\{1,\ldots,n\}$. Note that 
    \begin{align}
        \|\Delta_k & -hf_\delta(\theta_k,\tilde Y_{k-1})\|^2  = \bigl\|\Delta_{k-1}+h\bigl(f(\theta_k,\tilde Y_{k-1})-f(\theta_k,Y_{k-1})\bigr)\bigr\|^2 \notag \\ & = \|\Delta_{k-1}\|^2+2h\bigl\langle\Delta_{k-1},f(\theta_k,\tilde Y_{k-1})-f(\theta_k,Y_{k-1})\bigr\rangle+h^2\bigl\|f(\theta_k,\tilde Y_{k-1})-f(\theta_k,Y_{k-1})\bigr\|^2 \notag \\ & \leq \|\Delta_{k-1}\|^2 + C_1h\|\Delta_{k-1}\|^2 +L_1^2h^2 \|\Delta_{k-1}\|^{2\beta} \notag \\ & \leq (1+C_1h+L_1^2h^2)\|\Delta_{k-1}\|^2 +L_1^2h^2, \label{eq:noisy-lemma-1}
    \end{align}
    where in the passage to the third line, we used \eqref{itm:G3} and \eqref{itm:G4} with $U=\bar B_0(R^*)$; hence, 
    \begin{equation}
        C_1 = 2\max\{ 0, H_{\bar B_0(R^*)} \}
        \label{eq:C1-def}
    \end{equation}
    and $L_1$ is given by \eqref{eq:L1}. In the last passage, we used the inequality $x^{2\beta}\leq 1+x^2$ that holds for any $x\geq 0$ and $\beta \in (0,1)$.

    On the other hand,
    \begin{align}
        \|\Delta_k -hf_\delta(\theta_k,\tilde Y_{k-1})\|^2  & \geq \|\Delta_{k}\|^2 -2h\langle \Delta_k, f_\delta(\theta_k, \tilde Y_{k-1}) \rangle \notag \\ & \geq \|\Delta_{k}\|^2 -2h\|\Delta_k\|  \|f_\delta(\theta_k, \tilde Y_{k-1})\| \notag \\ & \geq \|\Delta_{k}\|^2 -h(\|\Delta_k\|^2+  \|f_\delta(\theta_k, \tilde Y_{k-1})\|^2) \notag \\ & \geq (1-h)\|\Delta_{k}\|^2 -h(1+R^*)^2\delta^2, \label{eq:noisy-lemma-2}
    \end{align}
    where in the second passage, we used the Cauchy-Schwarz inequality, in the third line -- the inequality $x^2+y^2\geq 2xy$ that holds for every $x,y\in\mathbb{R}$, and in the last line -- \eqref{eq:noise_function} and Corollary \ref{cor:R}.

    Combining \eqref{eq:noisy-lemma-1} with \eqref{eq:noisy-lemma-2}, we obtain 
    \begin{align}
       \|\Delta_{k}\|^2 & \leq \frac{1}{1-h}\bigl[ (1+C_1h+L_1^2h^2)\|\Delta_{k-1}\|^2 +L_1^2h^2 +h(1+R^*)^2\delta^2\bigr], \notag \\ & \leq (1+2h)(1+C_1h+L_1^2h^2)\|\Delta_{k-1}\|^2+2[L_1^2h^2 +h(1+R^*)^2\delta^2] \notag \\ & \leq (1+K_1h)\|\Delta_{k-1}\|^2+K_2h(h+\delta^2) \label{eq:noisy-lemma-3}
    \end{align}
    where the second and the third bounds follow from the assumption that $n>2(b-a)$ which implies that $0<h\leq \frac12$ and thus 
    $$0<\frac{1}{1-h}\leq 1+2h\leq 2.$$
    In the last passage, we substituted $K_1 = 2+2C_1+L_1^2$ and $K_2=2(\max\{L_1, 1+R^*\})^2$. By \eqref{eq:noisy-lemma-3} and the discrete Gronwall's lemma,
    \begin{align*}
       \|\Delta_{k}\|^2 & \leq (1+K_1h)^k\|\Delta_{0}\|^2+K_2h(h+\delta^2)\cdot \frac{(1+K_1h)^k-1}{K_1h} \\ & \leq e^{K_1(b-a)}\delta^2+\frac{K_2(e^{K_1(b-a)}-1)}{K_1}(h+\delta^2) \leq K_3(h+\delta^2), 
    \end{align*}
    where $K_3 = e^{K_1(b-a)}+\frac{K_2(e^{K_1(b-a)}-1)}{K_1}$. Thus,
    $$\|\Delta_{k}\|\leq \sqrt{K_3} \cdot \sqrt{h+\delta^2} \leq \sqrt{K_3} \bigl(h^{\frac12}+\delta\bigl),$$
    which completes the proof.
\end{proof}

The following theorem establishes the upper bound of the $L^p(\Omega)$-error of $y_n$, uniform over $n\in\mathbb{Z}_+$.

\begin{theorem}
    Let $p\in[2,\infty)$, $\eta\in\mathbb{R}^d$, and let the function $f\colon [a,b]\times\mathbb{R}^d\to\mathbb{R}^d$ satisfy the conditions \eqref{itm:G1}--\eqref{itm:G4}. Moreover, let $z, \tilde y_n\colon [a,b]\to\mathbb{R}^d$ be the exact solution to the IVP \eqref{eq:ivp} and the approximate solution given by \eqref{eq:barz}, respectively.
    Then there exists a positive constant $\tilde C$ such that for every $n\in\mathbb{Z}_+$ satisfying $n\geq 2(b-a)$ and every $\delta\in [0,1]$,
    \begin{equation}
        \Bigl\| \sup_{t\in [a,b]}\|\tilde y_n(t)-z(t)\| \Bigr\|_{L^p(\Omega)} \leq \tilde C\bigl( h^{\alpha+\frac12} + h^\beta+\delta + h^{\frac12}\mathbbm{1}_{(0,1]}(\delta)\mathbbm{1}_{(0,1)}(\beta) \bigr).
        \label{eq:bound_th1}
    \end{equation}
    \label{theorem}
\end{theorem}
    
\begin{proof}
    Note that
    \begin{equation}
        \sup_{t\in [a,b]}\|\tilde y_n(t)-z(t)\| \leq  \sup_{t\in [a,b]}\|z(t)-\bar z_n(t)\|+\sup_{t\in [a,b]}\|\bar z_n(t)-\tilde y_n(t)\|.
        \label{eq:step1}
    \end{equation}
    Since for each $k\in\{1,\ldots,n\}$, a function $\bar z_n-\tilde y_n$ is linear on $[t_{k-1},t_k]$, we get
    \begin{align*}
        \sup_{t\in [t_{k-1},t_k]}\|\bar z_n(t)-\tilde y_n(t)\| & = \max\bigl\{ \|\bar z_n(t_{k-1})-\tilde y_n(t_{k-1})\|,\|\bar z_n(t_{k})-\tilde y_n(t_{k})\| \bigr\} \\ & = \max\bigl\{\|z(t_{k-1})-\tilde Y_{k-1}\|, \|z(t_{k})-\tilde Y_{k}\|\bigr\}.
    \end{align*}    
    Thus,
    \begin{equation*}
        \sup_{t\in [a,b]}\|\bar z_n(t)-y_n(t)\| = \max_{k\in\{0,1,\ldots,n\}}\|z(t_{k})-\tilde Y_{k}\|\leq \max_{k\in\{0,1,\ldots,n\}}\|e_k\|+\max_{k\in\{0,1,\ldots,n\}}\|\tilde Y_k-Y_k\|.
    \end{equation*}
    Inserting this identity, inequality \eqref{eq:bound-z-lin-z} from Fact \ref{f:bound-z-lin-z}, and inequality \eqref{eq:max-Y-tilde-Y} from Lemma \ref{lemma:Y-tilde} into \eqref{eq:step1}, we obtain
    \begin{align*}
        \sup_{t\in [a,b]}\|\tilde y_n(t)-z(t)\| & \leq  L_1dh^{\alpha+1} +L_1dK^\beta(1+K_0)^\beta h^{\beta+1} +\tilde K \bigl( \delta + h^{\frac12}\mathbbm{1}_{(0,1]}(\delta)\mathbbm{1}_{(0,1)}(\beta) \bigr)\\ & \ \ \ \ \ \ \ \ \ \ \ \ \ \ \ \ \ \ \ \ \ \ \ \ \ \ \ \ \ \ \ \ \ \ \ \ \ \ \ \ \ \ \ \ \ \ \ \ \ \ \ \ \ \ \ \ \ \ \ \ \ \ \ \ \ \ \   +\max_{k\in\{0,1,\ldots,n\}}\|e_k\| \\ & \leq C_0\bigl( h^{\alpha+\frac12} + h^\beta +\delta + h^{\frac12}\mathbbm{1}_{(0,1]}(\delta)\mathbbm{1}_{(0,1)}(\beta)\bigr) +\max_{k\in\{0,1,\ldots,n\}}\|e_k\|,
    \end{align*}
    where $C_0 = \max\bigl\{\tilde K, L_1d(1+K_0)^\beta\max\{1,K^\beta\}\bigr\}$ and we used the fact that $h\leq 1$. As a result,
    \begin{equation}
        \Bigl\| \sup_{t\in [a,b]}\|\tilde y_n(t)-z(t)\| \Bigr\|_{L^p(\Omega)} \leq C_0( h^{\alpha+\frac12} + h^\beta +\delta + h^{\frac12}\mathbbm{1}_{(0,1]}(\delta)\mathbbm{1}_{(0,1)}(\beta) \bigr) + \Bigl\| \max_{k\in\{0,1,\ldots,n\}}\|e_k\| \Bigr\|_{L^p(\Omega)}.
        \label{eq:C0}
    \end{equation}
    With probability $1$, for every $k\in\{0,1,\ldots,n-1\}$ and $t\in [t_k,t_{k+1}]$,
    \begin{align*}
        \frac{\mathrm{d}}{\mathrm{d}t}\|z(t)-z_k(t)\|^2 & = 2 \langle z(t)-z_k(t),z'(t)-z'_k(t)\rangle = 2 \langle z(t)-z_k(t),f(t,z(t)) - f(t,z_k(t))\rangle \\ & \leq C_1 \|z(t)- z_k(t)\|^2.
    \end{align*}
    In the passage to the second line, we used \eqref{itm:G3} with $U=\bar B_0(R^*)$ and \eqref{eq:C1-def}. Since 
    \begin{equation*}
        z(t_k)-z_k(t_k) = z(t_k)-Y_k=e_k,
    \end{equation*}
    we obtain
    \begin{align*}
        \|z(t)-z_k(t)\|^2 & =  \|z(t_k)-z_k(t_k)\|^2+\int\limits_{t_k}^t\frac{\mathrm{d}}{\mathrm{d}s} \|z(s)- z_k(s)\|^2\,\mathrm{d}s \\ & \leq \|e_k\|^2 + C_1 \int\limits_{t_k}^t \|z(s)-z_k(s)\|^2\,\mathrm{d}s
    \end{align*}
    for $t\in [t_k,t_{k+1}]$. By Gronwall's inequality,
    \begin{equation}
        \|z(t)-z_k(t)\|^2 \leq \|e_k\|^2 e^{C_1 (t-t_k)}.
    \end{equation}
    As a result,
    \begin{equation}
        \|u_{k+1}\|^2 = \|z(t_{k+1})-z_k(t_{k+1})\|^2 \leq \|e_k\|^2 e^{C_1 (t_{k+1}-t_k)}= \|e_k\|^2 e^{C_1h}.
        \label{eq:u_bound}
    \end{equation}
    for every $k\in\{0,1,\ldots,n-1\}$. Furthermore, for $k\in\{0,1,\ldots,n-1\}$,
    \begin{align}
        e_{k+1} & = z(t_{k+1})-Y_{k+1} = z(t_{k+1}) - Y_k - hf(\theta_{k+1}, Y_k)= z(t_{k+1}) - z_k(t_k) - hf(\theta_{k+1}, Y_k) \notag \\ 
        & =  z(t_{k+1}) - \Bigl( z_k(t_{k+1}) - \int\limits_{t_k}^{t_{k+1}} f(s,z_k(s))\,\mathrm{d}s \Bigr) - hf(\theta_{k+1}, Y_k)  \notag \\ 
        & =  \underset{=u_{k+1}}{\underbrace{z(t_{k+1}) - z_k(t_{k+1})}} + \underset{=:A_k}{\underbrace{\int\limits_{t_k}^{t_{k+1}} \bigl(f(s,z_k(s))-f(s,Y_k)\bigr)\,\mathrm{d}s}} + \underset{=:B_k}{\underbrace{\int\limits_{t_k}^{t_{k+1}} f(s,Y_k)\,\mathrm{d}s-hf(\theta_{k+1}, Y_k)}}.  \label{eq:sum_uAB} 
    \end{align}
    Now we will establish upper bounds for $\|A_k\|$ and $\|B_k\|$. Fix $k\in\{0,1,\ldots,n-1\}$. Then
    \begin{equation}
        \|A_k\| \leq  \int\limits_{t_k}^{t_{k+1}} \bigl\|f(s,z_k(s))-f(s,Y_k)\bigr\|\,\mathrm{d}s \leq L_1 \int\limits_{t_k}^{t_{k+1}} \bigl\|z_k(s)-Y_k\bigr\|^\beta\,\mathrm{d}s,
        \label{eq:A_bound_1}
    \end{equation}
    where the last passage follows from \eqref{itm:G4} with $L_1 = L_{\bar B_0(R^*)}$; cf. \eqref{eq:L1}. By \eqref{itm:G2}, \eqref{eq:local_ivp}, and Corollary \ref{cor:R},
    \begin{align*}
        \bigl\|z_k(s)-Y_k\bigr\| & = \Bigl\| \int\limits_{t_k}^s f(u,z_k(u))\,\mathrm{d}u \Bigr\| \leq \int\limits_{t_k}^s \|f(u,z_k(u))\|\,\mathrm{d}u \leq \int\limits_{t_k}^sK(1+\|z_k(u)\|)\,\mathrm{d}u \\ & \leq Kh(1+R^*)
    \end{align*}
    for $s\in[t_k,t_{k+1}]$. Combining this inequality with \eqref{eq:A_bound_1}, we obtain
    \begin{equation}
        \|A_k\| \leq   L_1 \int\limits_{t_k}^{t_{k+1}} K^\beta h^\beta(1+R^*)^\beta\,\mathrm{d}s = L_1K^\beta (1+R^*)^\beta h^{\beta+1}= C_2  h^{\beta+1},
        \label{eq:A_bound_2}
    \end{equation}
    where $C_2 = L_1K^\beta (1+R^*)^\beta$.

    Note that $Y_k$ is $\mathcal{F}_k$-measurable, while $\theta_{k+1}$ is uniformly distributed over the interval $[t_k,t_{k+1}]$ and is independent of $\mathcal{F}_k$. Thus,
    \begin{equation}
        \mathbb{E}(B_k|\mathcal{F}_k) = \int\limits_{t_k}^{t_{k+1}} f(s,Y_k)\,\mathrm{d}s - h \mathbb{E}\bigl(f(\theta_{k+1}, Y_k)\bigr|Y_k\bigr) = 0.
        \label{eq:EB_k}
    \end{equation}
    Furthermore,
    \begin{align}
        \|B_k\| & =\Bigl\|\int\limits_{t_k}^{t_{k+1}} \bigl(f(s,Y_k)-f(\theta_{k+1},Y_k)\bigr)\,\mathrm{d}s\Bigr\| \leq \int\limits_{t_k}^{t_{k+1}} \bigl\|f(s,Y_k)-f(\theta_{k+1},Y_k)\bigr\|\,\mathrm{d}s \notag \\ &  \leq L_1 \int\limits_{t_k}^{t_{k+1}} |s-\theta_{k+1}|^{\alpha}\,\mathrm{d}s \leq L_1 h^{\alpha+1}, \label{eq:B_bound}
    \end{align}
    where the passage to the second line follows from \eqref{itm:G4}, and the last passage is obtained by the inequality $|s-\theta_{k+1}|\leq h$ that holds with probability $1$ for all $s\in[t_k,t_{k+1}]$.

    We go back to \eqref{eq:sum_uAB}:
    \begin{equation*}
        \|e_{k+1}\|^2 = \|u_{k+1}+A_k+B_k\|^2 = \|u_{k+1}+A_k\|^2+\|B_k\|^2 + 2\langle u_{k+1}+A_k, B_k \rangle.
    \end{equation*}
    Let $M_{k+1} = 2\langle u_{k+1}+A_k, B_k \rangle$ for $k\in\{0,1,\ldots,n-1\}$. For every $x,y\in\mathbb{R}^d$ and $h\in (0,1]$, we have, by applying the Young inequality, that
    \begin{align}
        \|x+y\|^2 & = \|x\|^2+2\langle x,y\rangle +\|y\|^2 \leq \|x\|^2+2\| x\| \|y\| +\|y\|^2\leq (1+h)\|x\|^2 + \Bigl(1+\frac1h\Bigr)\|y\|^2 \notag \\ & \leq (1+h)\|x\|^2+\frac2h \|y\|^2. \label{eq:x+y}
    \end{align}
    As a result,
    \begin{equation*}
        \|e_{k+1}\|^2\leq (1+h)\|u_{k+1}\|^2+\frac{2}{h}\|A_k\|^2+\|B_k\|^2+M_{k+1}.
    \end{equation*}
    The right-hand side can be further bounded by \eqref{eq:u_bound}, \eqref{eq:A_bound_2}, and \eqref{eq:B_bound}, leading to
    \begin{equation}
        \|e_{k+1}\|^2\leq (1+h)e^{C_1h}\|e_k\|^2 +2C_2^2 h^{2\beta+1}+L_1^2h^{2\alpha+2}+M_{k+1}.
        \label{eq:bound-ek}
    \end{equation}
    Let $\gamma:=(1+h)e^{C_1h}$. Note that 
    \begin{equation}
        \mathbb{E}(M_{k+1}|\mathcal{F}_k) = \langle u_{k+1}+A_k, \mathbb{E}(B_k|\mathcal{F}_k) \rangle = 0
        \label{eq:EM_k}
    \end{equation}
    by  $\mathcal{F}_k$-measurability of $u_{k+1}+A_k$ and by \eqref{eq:EB_k}. As a result, by \eqref{eq:bound-ek}, 
    \begin{equation*}
        \mathbb{E}\|e_{k+1}\|^2 \leq \gamma\mathbb{E}\|e_k\|^2 +2C_2^2h^{2\beta+1}+L_1^2 h^{2\alpha+2}.
    \end{equation*}
    By discrete Gronwall's lemma,
    \begin{equation*}
        \mathbb{E}\|e_{k}\|^2 \leq \mathbb{E}\|e_0\|^2 \gamma^k+\frac{\gamma^k-1}{\gamma-1}\bigl(2C_2^2 h^{2\beta+1}+L_1^2 h^{2\alpha+2}\bigr).
    \end{equation*}
     Since $e_0=0$, $\gamma-1=(1+h)e^{C_1h}-1\geq (1+h)-1=h$, and 
     \begin{equation}
        1\leq \gamma^k=(1+h)^ke^{C_1hk}\leq \Bigl(1+\frac{b-a}{n}\Bigr)^n e^{C_1(b-a)} \leq e^{(C_1+1)(b-a)} =: C_3,
        \label{eq:C3}
    \end{equation}
    we obtain
    \begin{equation}
        \mathbb{E}\|e_{k}\|^2 \leq \frac{C_3-1}{h}\bigl(2C_2^2 h^{2\beta+1}+L_1^2 h^{2\alpha+2}\bigr) \leq C_4 \bigl( h^{2\alpha+1} + h^{2\beta} \bigr)
        \label{eq:Eek}
    \end{equation}
    for $k\in\{0,1,\ldots,n\}$, where $C_4: = (C_3-1)\cdot\max\{2C_2^2, L_1^2\}$.
    
    On the other hand, let $E_k := \gamma^{-k}\|e_k\|^2$ for $k\in\{0,1,\ldots,n\}$. Then \eqref{eq:bound-ek} can be restated as
    \begin{equation*}
        E_{k+1}\leq E_k +(2C_2^2 h^{2\beta+1}+L_1^2h^{2\alpha+2}+M_{k+1})\gamma^{-(k+1)}
    \end{equation*}
    for every $k\in\{0,1,\ldots,n-1\}$. Hence,
    \begin{align}
        E_k & = E_0 + \sum_{j=0}^{k-1} (E_{j+1}-E_j) \leq \sum_{j=0}^{k-1} (2C_2^2 h^{2\beta+1}+L_1^2h^{2\alpha+2}+M_{j+1})\gamma^{-(j+1)} \notag \\ & \leq \sum_{j=0}^{k-1} (2C_2^2 h^{2\beta+1}+L_1^2h^{2\alpha+2}) + \sum_{j=0}^{k-1} M_{j+1}\gamma^{-(j+1)} \notag \\ & \leq C_5h^{2\beta} + L_1^2(b-a)h^{2\alpha+1} + N_k, \label{eq:Ek}
    \end{align}
    where the passage to the second line follows from the fact that $\gamma\geq 1$ and thus $\gamma^{-(j+1)}\in (0,1]$ for all $j\in\{0,1,\ldots,n-1\}$. In the third line, $C_5:= 2C_2^2(b-a)$ and $\displaystyle N_k := \sum_{j=0}^{k-1} M_{j+1}\gamma^{-(j+1)}$. 
    Moreover, since $\gamma^{k}E_k = \|e_k\|^2$ and by \eqref{eq:C3}, we have
    \begin{equation*}
        C_3\cdot \max_{k\in\{0,1,\ldots,n\}} E_k \geq \max_{k\in\{0,1,\ldots,n\}}\|e_k\|^2.
    \end{equation*}
    Combining this inequality with \eqref{eq:Ek},
    \begin{align}
        \max_{k\in\{0,1,\ldots,n\}}\|e_k\|^2  & \leq C_3C_5h^{2\beta} + C_3L_1^2(b-a)h^{2\alpha+1} + C_3\max_{k\in\{1,\ldots,n\}} |N_k| \notag \\ & \leq C_6(h^{2\beta} + h^{2\alpha+1}) + C_3\max_{k\in\{1,\ldots,n\}} |N_k|, \label{eq:max-bound-as}
    \end{align}
    where $C_6:=C_3\cdot \max\{C_5, L_1^2(b-a)\}$.

    Let $q:=\frac{p}{2}\geq 1$. By \eqref{eq:max-bound-as},
    \begin{align}
        \Bigl\| \max_{k\in\{0,1,\ldots,n\}}\|e_k\| \Bigr\|_{L^p(\Omega)}^2 & = \Bigl[\mathbb{E}\Bigl(\max_{k\in\{0,1,\ldots,n\}}\|e_k\|\Bigr)^p\Bigr]^{\frac{2}{p}} \notag\\ & = \Bigl[\mathbb{E}\Bigl(\max_{k\in\{0,1,\ldots,n\}}\|e_k\|^2\Bigr)^q\Bigr]^{\frac{1}{q}} \notag \\ & = \Bigl\|  \max_{k\in\{0,1,\ldots,n\}}\|e_k\|^2\Bigr\|_{L^q(\Omega)} \notag \\ & \leq \Bigl\| C_6(h^{2\beta} + h^{2\alpha+1}) + C_3\max_{k\in\{1,\ldots,n\}} |N_k| \Bigr\|_{L^q(\Omega)} \notag \\ & \leq C_6(h^{2\beta} + h^{2\alpha+1}) + C_3\Bigl\| \max_{k\in\{1,\ldots,n\}} |N_k| \Bigr\|_{L^q(\Omega)}, \label{eq:max-case2}
    \end{align}
    Recall that $N_k$ was defined under \eqref{eq:bound-ek} as $\displaystyle N_k = \sum_{j=0}^{k-1} M_{j+1}\gamma^{-(j+1)}$. Let additionally $N_0=0$. Then $(N_k)_{k=0}^n$ is an $(\mathcal{F}_k)_{k=0}^n$-martingale. Indeed, for each $k\in\{1,\ldots,n-1\}$,
    \begin{equation}
        \mathbb{E}(N_{k+1}|\mathcal{F}_k) = \sum_{j=0}^k \mathbb{E}(M_{j+1}|\mathcal{F}_k) \gamma^{-(j+1)} = \sum_{j=0}^{k-1} M_{j+1}\gamma^{-(j+1)} + \mathbb{E}(M_{k+1}|\mathcal{F}_k)\gamma^{-(k+1)}=N_k,
        \label{eq:N-martingale}
    \end{equation}
where we used the fact that $M_j$ is $\mathcal{F}_k$-measurable for $j\leq k$, and applied \eqref{eq:EM_k}. The quadratic variation of this martingale is
    \begin{equation}
        [N]_k = |N_0|^2 + \sum_{j=1}^k |N_{j}-N_{j-1}|^2 = \sum_{j=1}^k M_j^2 \gamma^{-2j}.
        \label{eq:quadratic_variation}
    \end{equation}
    By Burkholder's inequality, there exists a positive constant $C_q$ such that for each $n\in\mathbb{Z}_+$,
    \begin{equation}
        \Bigl\| \max_{k\in\{1,\ldots,n\}} |N_k| \Bigr\|_{L^q(\Omega)} \leq C_q \bigl\| [N]^{\frac{1}{2}}_n \bigr\|_{L^q(\Omega)}.
        \label{eq:burkholder}
    \end{equation}
    Let us bound $[N]_n$ given by \eqref{eq:quadratic_variation}:
    \begin{align}
        [N]_n & = \sum_{k=1}^n M_k^2 \gamma^{-2k} \leq \sum_{k=1}^n M_k^2 = \sum_{k=0}^{n-1} 4\langle u_{k+1}+A_{k}, B_k \rangle^2 \leq 4\sum_{k=0}^{n-1} \|u_{k+1}+A_{k}\|^2\|B_k\|^2 \notag \\ & \leq 4L_1^2h^{2\alpha+2}\sum_{k=0}^{n-1} \|u_{k+1}+A_{k}\|^2,
        \label{eq:varN-bound}
    \end{align}
    where the first inequality is due to the fact that $\gamma^{-2k}\leq 1$ for each $k\in\{1,\ldots,n\}$, cf. \eqref{eq:C3}, and the second inequality is an application of the Cauchy-Schwarz inequality. The passage to the second line is due to \eqref{eq:B_bound}.
    Note that
    \begin{align*}
        \|u_{k+1}+A_{k}\|^2 & \leq \bigl(\|u_{k+1}\|+\|A_{k}\|\bigr)^2 \leq 2\bigl(\|u_{k+1}\|^2+\|A_{k}\|^2\bigr) \notag \\ & \leq 2e^{C_1h}\|e_k\|^2+2C_2^2h^{2\beta+2},
    \end{align*}
    where in the passage to the second line, we used \eqref{eq:u_bound} and \eqref{eq:A_bound_2}. As a result,
    \begin{align*}
        \Bigl( \sum_{k=0}^{n-1} \|u_{k+1}+A_{k}\|^2 \Bigr)^{\frac{1}{2}} & \leq \Bigl( 2e^{C_1h}\sum_{k=0}^{n-1} \|e_k\|^2+n\cdot 2C_2^2h^{2\beta+2} \Bigr)^{\frac{1}{2}} \\ & \leq \Bigl( 2e^{C_1(b-a)}\cdot n \max_{k\in\{0,1,\ldots,n\}}\|e_k\|^2+n\cdot 2C_2^2h^{2\beta+2} \Bigr)^{\frac{1}{2}} \\ & \leq \Bigl( 2e^{C_1(b-a)}\cdot n \max_{k\in\{0,1,\ldots,n\}}\|e_k\|^2\Bigr)^{\frac{1}{2}}+\Bigl(2C_2^2(b-a)h^{2\beta+1} \Bigr)^{\frac{1}{2}} \\ & = \sqrt{2e^{C_1(b-a)}}\cdot\sqrt{n}\cdot\max_{k\in\{0,1,\ldots,n\}}\|e_k\|+\sqrt{2(b-a)}C_2h^{\beta+\frac12},
    \end{align*}
    where in the passage to the third line, we used an inequality $\sqrt{x+y}\leq\sqrt{x}+\sqrt{y}$, which holds for any $x,y\geq 0$.    
    
    Combining this bound with \eqref{eq:varN-bound}, we obtain
    \begin{align*}
        \bigl\| [N]^{\frac{1}{2}}_n \bigr\|_{L^q(\Omega)} & \leq 2L_1h^{\alpha+1}\Bigl\| \Bigl( \sum_{k=0}^{n-1} \|u_{k+1}+A_{k}\|^2 \Bigr)^{\frac{1}{2}} \Bigr\|_{L^q(\Omega)} \\ & \leq 2L_1h^{\alpha+1} \Bigl( \sqrt{2e^{C_1(b-a)}}\cdot\sqrt{n}\cdot \Bigl\| \max_{k\in\{0,1,\ldots,n\}}\|e_k\|\Bigr\|_{L^q(\Omega)} +\sqrt{2(b-a)}C_2h^{\beta+\frac12}\Bigr)  \\ & \leq C_7\Bigl(h^{\alpha+\frac12}\cdot \Bigl\| \max_{k\in\{0,1,\ldots,n\}}\|e_k\|\Bigr\|_{L^q(\Omega)} +h^{\alpha+\beta+\frac32} \Bigr) \\ & \leq C_7\Bigl(h^{\alpha+\frac12}\cdot \Bigl\| \max_{k\in\{0,1,\ldots,n\}}\|e_k\|\Bigr\|_{L^p(\Omega)} +h^{\alpha+\beta+\frac32}\Bigr),
    \end{align*}
    where
    $$C_7:= 2L_1\sqrt{2(b-a)}\max\bigl\{\sqrt{e^{C_1(b-a)}}, C_2\bigr\}$$
    and in the last passage, we used the fact that $p>q$.
    Let us insert this bound into \eqref{eq:burkholder} and then into \eqref{eq:max-case2}. We obtain:
    \begin{align*}
        \Bigl\| \max_{k\in\{0,1,\ldots,n\}}\|e_k\| \Bigr\|^2_{L^p(\Omega)} & \leq C_6(h^{2\beta} + h^{2\alpha+1}) + C_3C_q \bigl\| [N]^{\frac{1}{2}}_n \bigr\|_{L^q(\Omega)} \\ & \leq C_6(h^{2\beta} + h^{2\alpha+1}) + C_3C_7C_q \Bigl( h^{\alpha+\frac12}\cdot \Bigl\| \max_{k\in\{0,1,\ldots,n\}}\|e_k\|\Bigr\|_{L^p(\Omega)} +h^{\alpha+\beta+\frac32} \Bigr).
    \end{align*}
    Since $h\in (0,1]$ and $h^{\alpha+\beta+\frac32}\leq h^{\alpha+\beta+\frac12}\leq \frac12\bigl(h^{2\alpha+1}+h^{2\beta}\bigr)$,
    \begin{align*}
        \Bigl\| \max_{k\in\{0,1,\ldots,n\}}\|e_k\| \Bigr\|^2_{L^p(\Omega)} & \leq C_8(h^{2\beta} + h^{2\alpha+1}) + C_8  h^{\alpha+\frac12}\cdot \Bigl\| \max_{k\in\{0,1,\ldots,n\}}\|e_k\|\Bigr\|_{L^p(\Omega)} \\ & \leq C_8^2(h^{2\beta} + h^{2\alpha+1}) + C_8  h^{\alpha+\frac12}\cdot \Bigl\| \max_{k\in\{0,1,\ldots,n\}}\|e_k\|\Bigr\|_{L^p(\Omega)},
    \end{align*}
    where
    $$C_8:=\max\Bigl\{ C_6+\frac12C_3C_7C_q, C_3C_7C_q, 1\Bigr\}.$$
    Hence,
    \begin{align*}
        \frac54 \Bigl\| \max_{k\in\{0,1,\ldots,n\}}\|e_k\| \Bigr\|^2_{L^p(\Omega)} & \leq \Bigl( C_8h^{\alpha+\frac12} + \frac12 \Bigl\| \max_{k\in\{0,1,\ldots,n\}}\|e_k\| \Bigr\|_{L^p(\Omega)} \Bigr)^2 + C_8^2h^{2\beta} \\ & \leq \Bigl( C_8h^{\alpha+\frac12} + \frac12 \Bigl\| \max_{k\in\{0,1,\ldots,n\}}\|e_k\| \Bigr\|_{L^p(\Omega)} + C_8h^\beta\Bigr)^2.
    \end{align*}
    As a result,
    \begin{equation*}
        \frac{\sqrt{5}-1}{2} \Bigl\| \max_{k\in\{0,1,\ldots,n\}}\|e_k\| \Bigr\|_{L^p(\Omega)} \leq C_8\bigl(h^{\alpha+\frac12}+h^\beta\bigr)
    \end{equation*}
    and finally, inserting the bound obtained above into \eqref{eq:C0},
    \begin{equation*}
        \Bigl\| \sup_{t\in [a,b]}\|\tilde y_n(t)-z(t)\| \Bigr\|_{L^p(\Omega)} \leq \tilde C\bigl(h^{\alpha+\frac12}+h^\beta+\delta + h^{\frac12}\mathbbm{1}_{(0,1]}(\delta)\mathbbm{1}_{(0,1)}(\beta)\bigr)
    \end{equation*}
    for $\displaystyle \tilde C = C_0+\frac{2}{\sqrt{5}-1} C_8$.
\end{proof}

\begin{remark}
 The order of convergence for the randomized Euler scheme is:
 \begin{itemize}
     \item $O(h^{\alpha+\frac12}+h^\beta)$ if $\delta=0$ (i.e., under exact information);
     \item $O(h^{\alpha+\frac12}+h^\beta+\delta)$ if $\delta \in (0,1]$ (i.e., under noisy information) and $\beta \in (0,\frac12)$ or $\beta=1$;
     \item $O(h^{\alpha+\frac12}+h^{\frac12}+\delta)$ if $\delta \in (0,1]$ (i.e., under noisy information) and $\beta \in [\frac12,1)$.
 \end{itemize}
\end{remark}

\begin{remark}
 For the deterministic explicit Euler scheme, the error bound is of order $O(h^{\alpha} + h^\beta)$. Thus, randomization improves the convergence rate of the error component associated with the H\"older regularity in time from $\alpha$ to $\alpha+\frac12$.
\end{remark}

\begin{remark}
 Randomization introduced in the Euler scheme applies exclusively to the time variable. As a result, for autonomous ordinary differential equations, the randomized and deterministic Euler schemes are equivalent, implying an error of order $O(h^\beta)$ under exact information.
\end{remark}

\section{Numerical experiments} \label{sec:experiments}

\begin{figure}[!ht]
    \centering
    \includegraphics[width=0.7\linewidth]{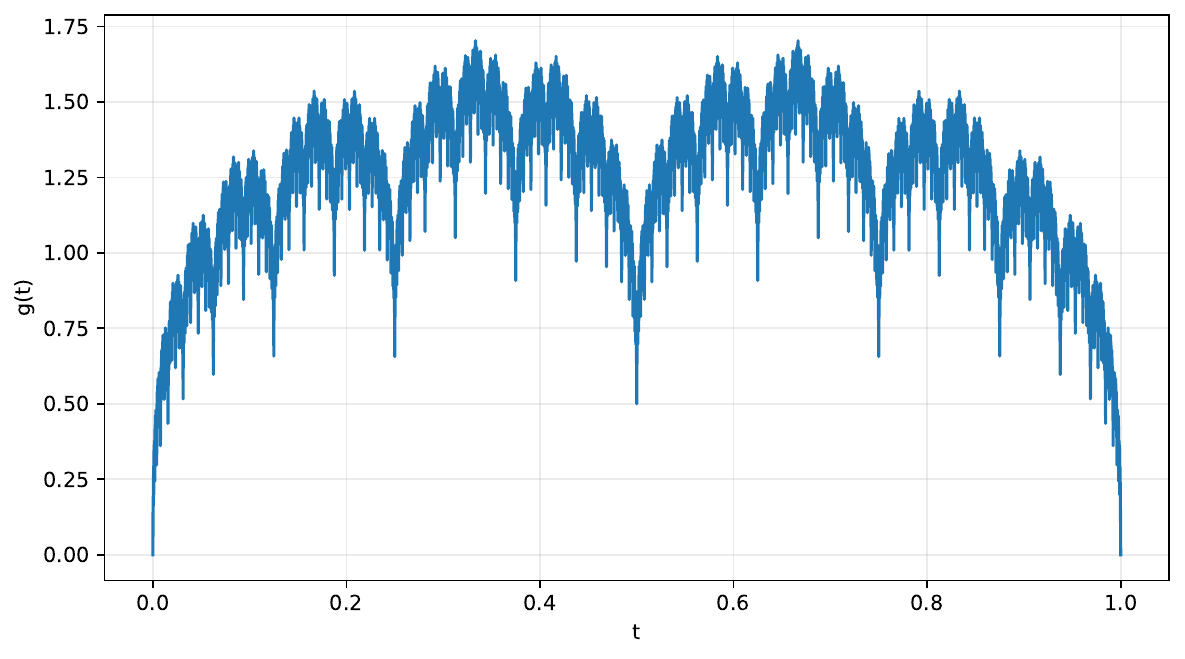}
    \caption{Takagi-Landsberg function $g_{\alpha,M}$ for $\alpha = 0.3$ and $M=14$.}
    \label{fig:tl}
\end{figure}

In this section, we compare the performance of the randomized Euler scheme with that of its deterministic counterpart. We consider initial value problems (IVPs) whose right-hand sides involve Takagi–Landsberg functions defined by
$$\phi(t)=\min(\{t\},1-\{t\}),\qquad g_{\alpha,M}(t)=\sum_{m=0}^{M}2^{-\alpha m}\phi(2^m t),$$
where $M\in\mathbb{Z}_+$ and $\alpha\in (0,1)$ are parameters, and $\{t\}$ denotes the fractional part of $t\in \mathbb{R}$. Takagi-Landsberg functions are H\"older continuous with exponent $\alpha$; see Remark 2.9 in \cite{TLfunctions}. Figure \ref{fig:tl} illustrates an example of a Takagi-Landsberg function.

\begin{figure}[h!]
    \centering
    \includegraphics[width=0.7\linewidth]{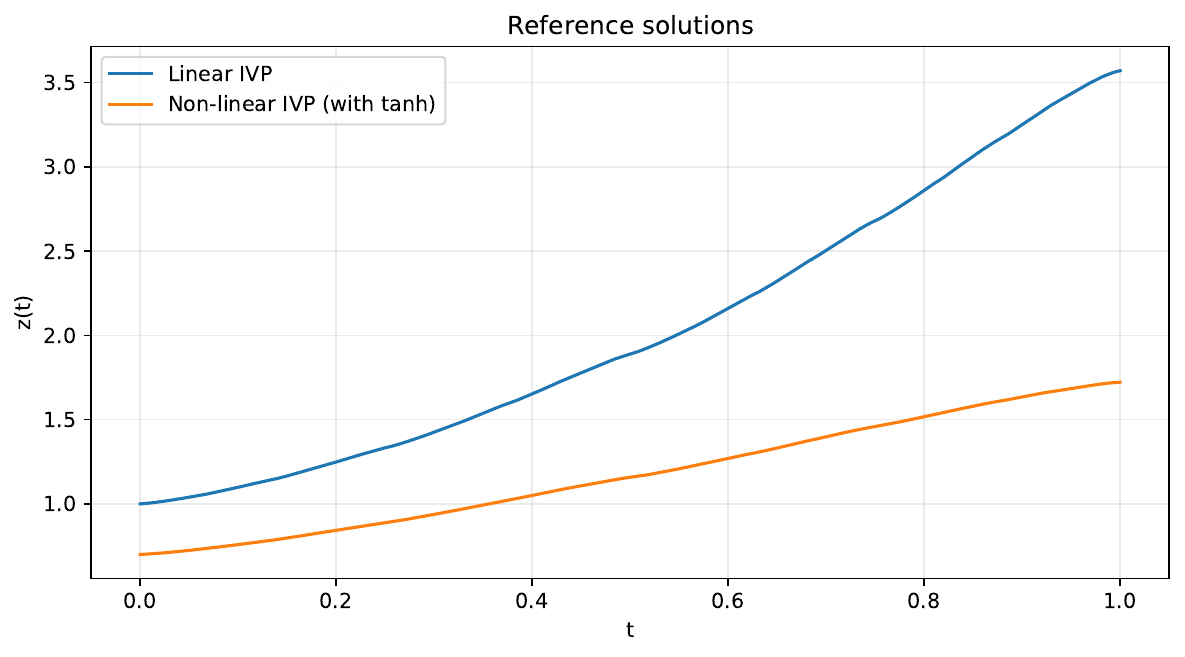}
    \caption{Exact solutions (up to the numerical integration of the Takagi–Landsberg function) to problems \eqref{eq:ivp1} and \eqref{eq:ivp2} for $\alpha=0.3$.}
    \label{fig:ref_sol}
\end{figure}

\noindent We consider three test IVPs:
\begin{enumerate}
    \item linear IVP
    \begin{equation}
        \left\{ 
        \begin{array}{l}
             z'(t) = g_{\alpha}(t)\cdot z(t), \ \ t\in [0,1],\\
             z(0) = 1, 
        \end{array}\right.
        \label{eq:ivp1}
    \end{equation}
    \item nonlinear IVP with $\tanh$:
    \begin{equation}
        \left\{ 
        \begin{array}{l}
             z'(t) = g_{\alpha}(t)\cdot \tanh z(t), \ \ t\in [0,1],\\
             z(0) = 0.7, 
        \end{array}\right.
        \label{eq:ivp2}
    \end{equation}
    \item and nonlinear IVP with $\cos$:
    \begin{equation}
        \left\{ 
        \begin{array}{l}
             z'(t)=-0.3\,g_{\alpha}(t)\operatorname{sgn}(z(t))|z(t)|^{2/3}|\cos z(t)|, \ \ t\in [-3.5,3.5],\\
             z(-3.5) = 3, 
        \end{array}\right.
        \label{eq:ivp3}
    \end{equation}
\end{enumerate}
where $g_\alpha = g_{\alpha, 14}$ is a Takagi-Landsberg function for $M=14$. The exact solution for problem \eqref{eq:ivp1} is given as
$$z(t) = \exp (G(t)), \quad G(t) = \int\limits_0^t g_\alpha(u)\,\mathrm{d}u,$$
whereas for problem \eqref{eq:ivp2}:
$$z(t) = \text{arsinh}\bigl( \sinh z(0)\cdot \exp (G(t)) \bigr).$$
These solutions are shown for $\alpha=0.3$ in Figure \ref{fig:ref_sol}.

\begin{figure}[!ht]
    \centering
    \includegraphics[width=0.7\linewidth]{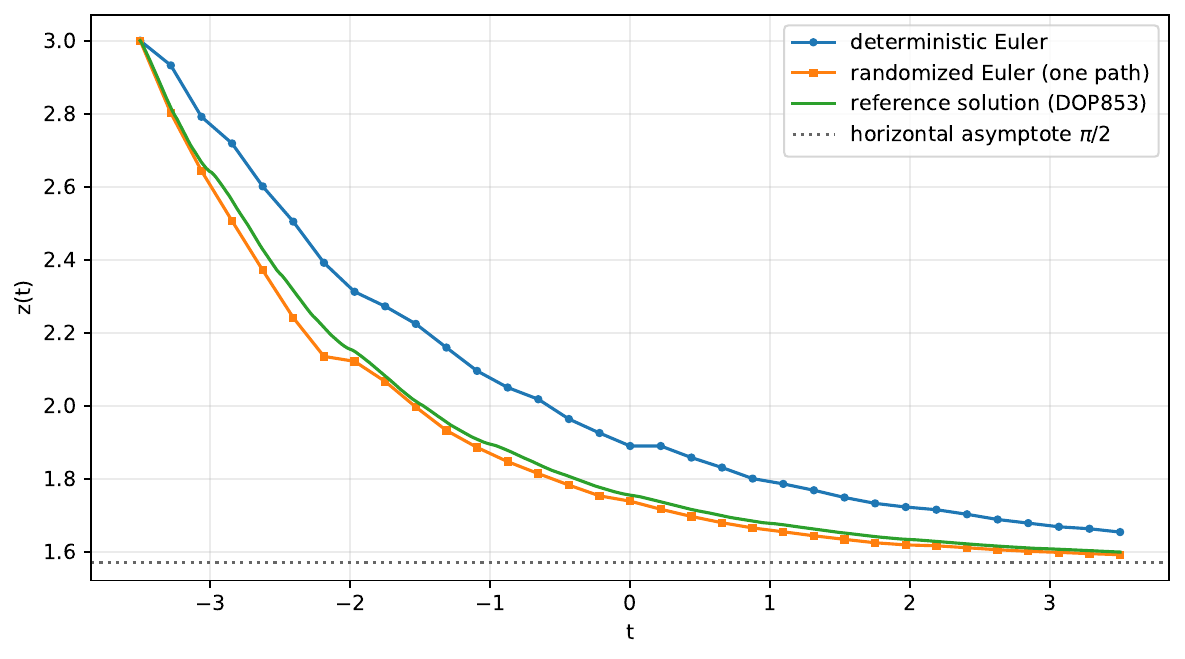}
    \caption{Reference solution to problem \eqref{eq:ivp3} for $\alpha=0.3$.}
    \label{fig:ref_sol2}
\end{figure}

For problem \eqref{eq:ivp3}, we first solve an autonomous IVP $v'(s)=F(v(s))$, $v(0)=3$, and then take $z(t)=v(G(t))$, where 
$$F(x)=-0.3\operatorname{sgn}(x)|x|^{2/3}|\cos x|,\qquad G_\alpha(t)=\int_{-3.5}^t g_{\alpha}(s)\,\mathrm{d}s.$$
Note that
$$z'(t) = \bigl( v(G_\alpha(t)) \bigr)' = v'(G_\alpha(t))\cdot G_\alpha'(t) = F(v(G_\alpha(t)))\cdot g_\alpha(t) = F(z(t))\cdot g_\alpha(t).$$
Thus, $z$ is an exact solution to \eqref{eq:ivp3}. To solve $v'(s)=F(v(s))$ numerically, we use SciPy's DOP853 solver. We integrate the function $G_a$ on a grid that includes all points of non-differentiability. The resulting reference solution is shown in Figure \ref{fig:ref_sol2}, together with approximate solutions obtained using the deterministic and randomized Euler schemes; for the latter, a single sample trajectory is shown. The solution approaches the horizontal asymptote $y=\pi/2$.

\begin{figure}[!ht]
    \centering
    \includegraphics[width=0.7\linewidth]{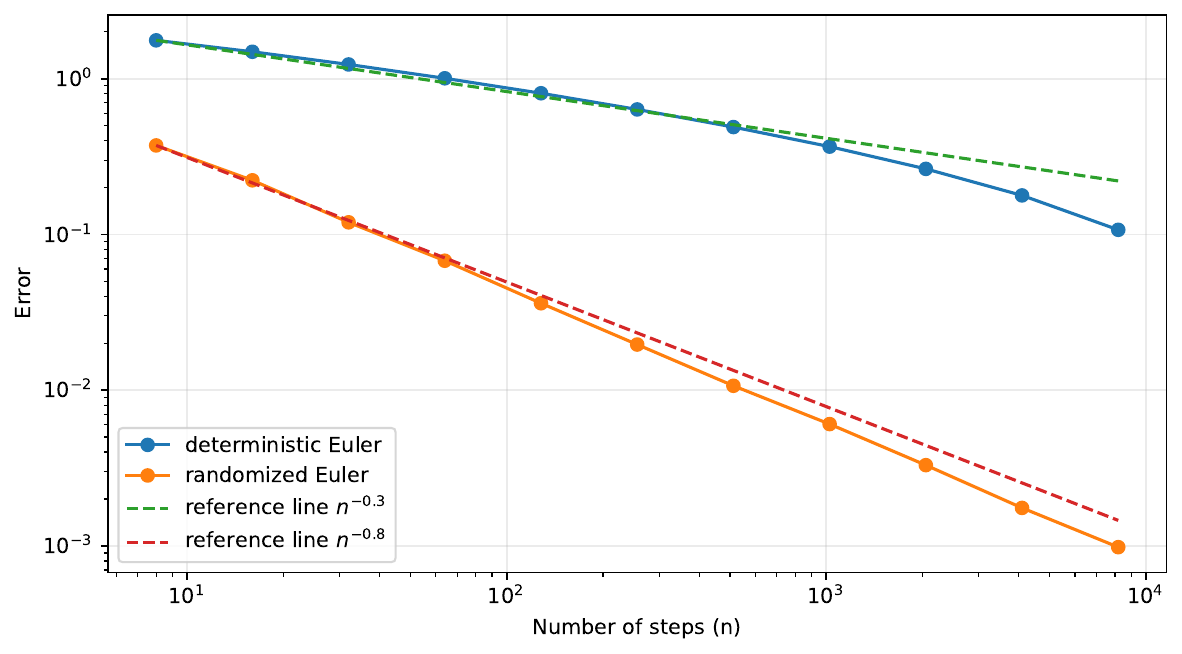}
    \caption{Errors of the deterministic and randomized Euler schemes applied to the linear problem \eqref{eq:ivp1} for $\alpha=0.3$.}
    \label{fig:linear}
\end{figure}

\begin{figure}[!ht]
    \centering
    \includegraphics[width=0.7\linewidth]{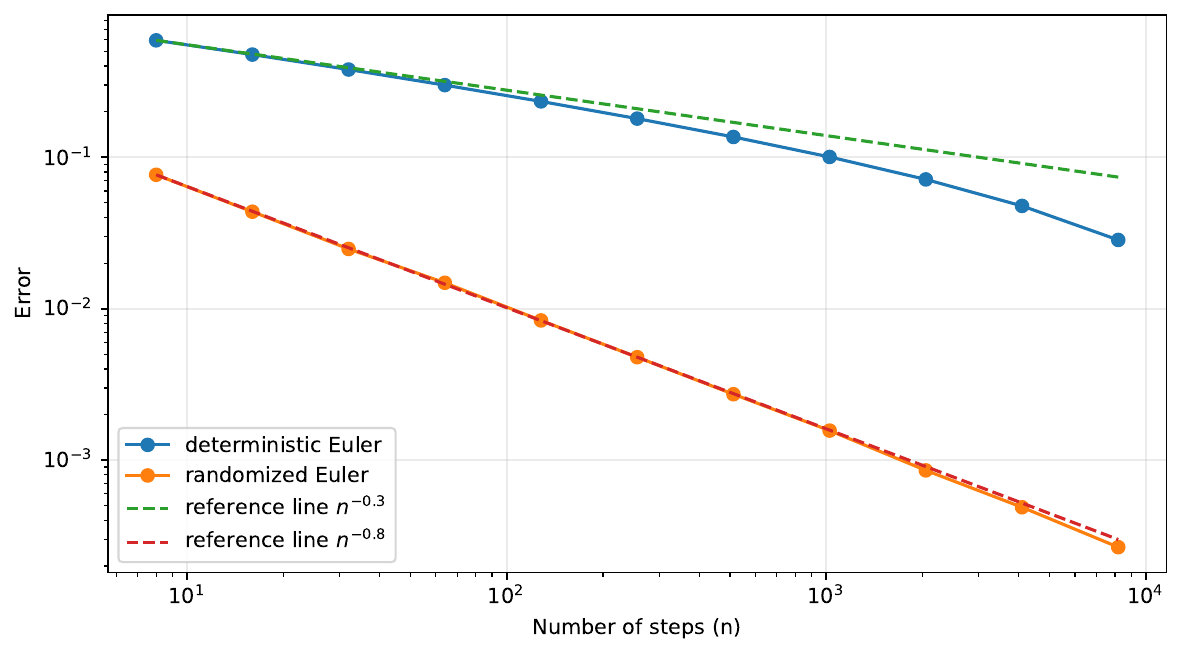}
    \caption{Errors of the deterministic and randomized Euler schemes applied to the nonlinear problem \eqref{eq:ivp2} for $\alpha=0.3$.}
    \label{fig:nonlinear}
\end{figure}

\begin{figure}[!ht]
    \centering
    \includegraphics[width=0.7\linewidth]{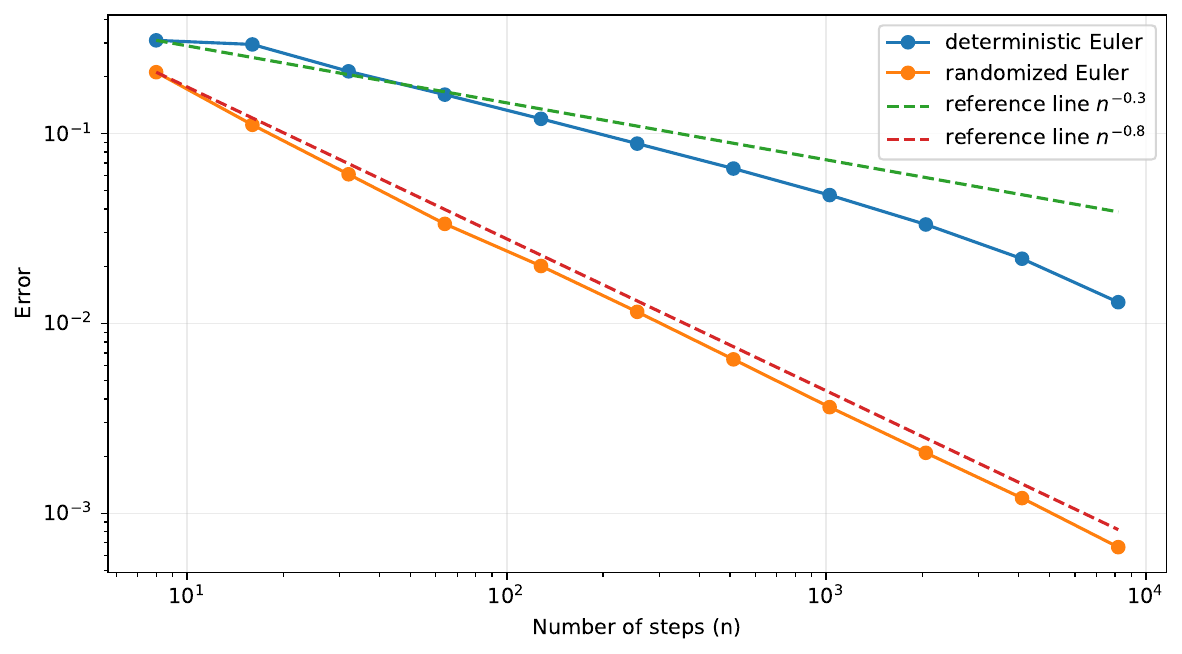}
    \caption{Errors of the deterministic and randomized Euler schemes applied to the nonlinear problem \eqref{eq:ivp3} for $\alpha=0.3$.}
    \label{fig:nonlinear2}
\end{figure}

Figures \ref{fig:linear}, \ref{fig:nonlinear}, and \ref{fig:nonlinear2} show the errors of the randomized and deterministic Euler schemes as functions of the number of time steps, $n\in \{2^3, 2^4,\ldots, 2^{13}\}$, for the linear problem \eqref{eq:ivp1} and the nonlinear problems \eqref{eq:ivp2} and \eqref{eq:ivp3}, respectively. In each of the problems, we set $\alpha=0.3$. Axes are displayed on a logarithmic scale. 

For the randomized scheme, the error is estimated using $1,\!000$ Monte Carlo simulations. The reference lines $n^{-\alpha}=n^{-0.3}$ and $n^{-\alpha-1/2}=n^{-0.8}$ correspond to the theoretical convergence rates of the deterministic and randomized Euler schemes, respectively (see Theorem \ref{theorem} and Theorem 3.1 in \cite{CzyzMorPrzyb}).

\begin{table}[!ht]
    \centering
    \resizebox{\textwidth}{!}{
    \begin{tabular}{c|c|c|c|c|c|c|c|c|c}
    & \multicolumn{3}{c|}{linear IVP \eqref{eq:ivp1}} & \multicolumn{3}{c|}{nonlinear IVP \eqref{eq:ivp2}} & \multicolumn{3}{c}{nonlinear IVP \eqref{eq:ivp3}} \\ \cline{2-10}
    & $\alpha=0.2$ & $\alpha=0.3$ & $\alpha=0.4$ & $\alpha=0.2$ & $\alpha=0.3$ & $\alpha=0.4$ & $\alpha=0.2$ & $\alpha=0.3$ & $\alpha=0.4$ \\ \hline
    deterministic Euler & $0.297$ & $0.389$ & $0.480$ & $0.351$ & $0.421$ & $0.499$ & $0.409$ & $0.457$ & $0.524$ \\ \hline 
    randomized Euler & $0.808$ & $0.864$ & $0.926$ & $0.728$ & $0.813$ & $0.901$ & $0.764$ & $0.821$ & $0.887$
    \end{tabular}}
    \caption{Empirical convergence rates of the deterministic and randomized Euler schemes applied to problems \eqref{eq:ivp1} and \eqref{eq:ivp2} for $\alpha\in\{0.2, 0.3, 0.4\}$. Rates were estimated from the errors for $n\in \{2^3, 2^4,\ldots, 2^{13}\}$.}
    \label{tab:rates}
\end{table}

The empirical convergence rates, estimated over the range $n\in[8,8192]$, are close to the theoretical predictions: the deterministic scheme exhibits a rate consistent with the H\"older exponent $\alpha$, whereas the randomized scheme achieves an improvement of $1/2$ in the rate. This observation is further supported by the results reported in Table \ref{tab:rates}, where different values of $\alpha$ are tested. For problem \eqref{eq:ivp3}, however, we have $\beta=\frac23$. Consequently, the empirical convergence rate of the randomized scheme exceeds the theoretical bound $\min\{\alpha+\frac12,\beta\}=\frac23$. This can be explained by observing that, although the right-hand side of \eqref{eq:ivp3} is only $\frac23$-H\"older continuous with respect to the space variable, its regularity is apparently sufficient to prevent the convergence rate from deteriorating to the minimal level established in Theorem \ref{theorem}.

\begin{figure}[!ht]
    \centering
    \includegraphics[width=0.7\linewidth]{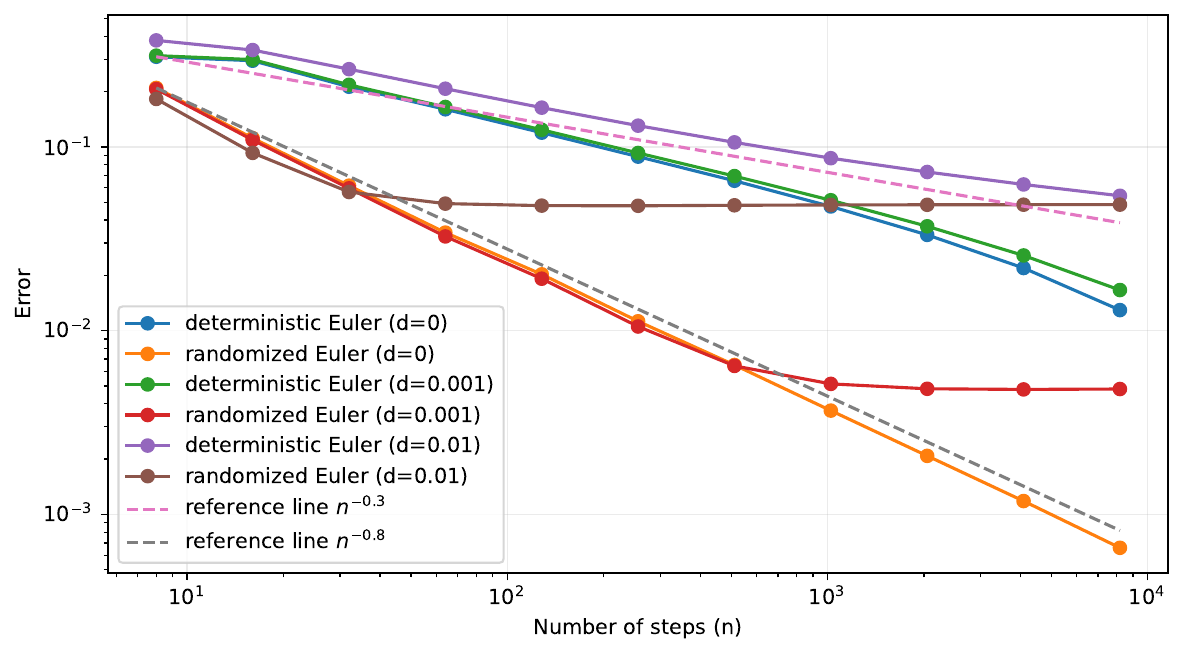}
    \caption{Errors of the deterministic and randomized Euler schemes applied to the nonlinear problem \eqref{eq:ivp3} for $\alpha=0.3$ with additive noise and different levels of the precision parameter $\delta\in\{0,10^{-3},10^{-2}\}$.}
    \label{fig:noisy}
\end{figure}

\begin{figure}[!ht]
    \centering
    \includegraphics[width=0.7\linewidth]{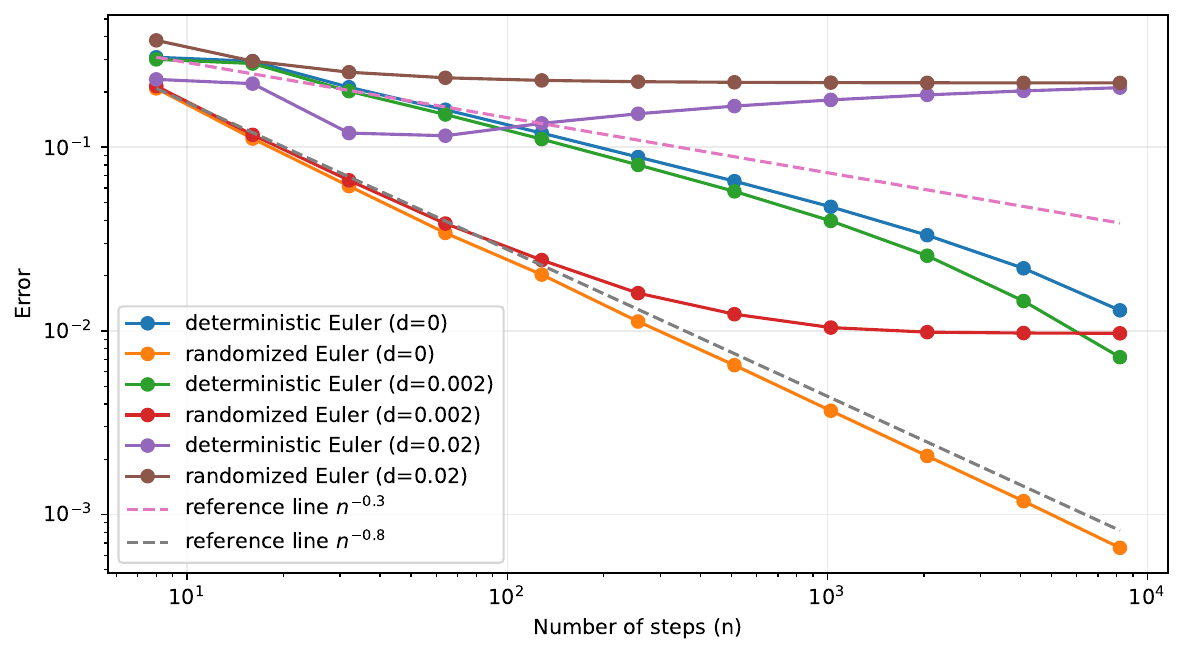}
    \caption{Errors of the deterministic and randomized Euler schemes applied to the nonlinear problem \eqref{eq:ivp3} for $\alpha=0.3$ with subtractive noise and different levels of the precision parameter $\delta\in\{0,10^{-3},10^{-2}\}$.}
    \label{fig:noisy_negative}
\end{figure}

Finally, Figures \ref{fig:noisy} and \ref{fig:noisy_negative} illustrate the impact of the precision parameter for problem \eqref{eq:ivp3}. We consider inexact information with two levels of the precision parameter: $\delta=0.002$ and $\delta=0.02$. We simulate the noise in two ways:
\begin{itemize}
    \item \textit{additive noise}, where, at each iteration of \eqref{eq:noisy_rand_ee}, we set $$\tilde Y_j = \tilde Y_{j-1} + h\cdot\bigl[f(\theta_j, \tilde Y_{j-1}) + \delta\cdot(1+|\tilde Y_{j-1}|) \bigr];$$
    \item \textit{subtractive noise}, where, at each iteration of \eqref{eq:noisy_rand_ee}, we set $$\tilde Y_j = \tilde Y_{j-1} + h\cdot\bigl[f(\theta_j, \tilde Y_{j-1}) - \delta\cdot(1+|\tilde Y_{j-1}|) \bigr].$$
\end{itemize}
The results for additive and subtractive noise are presented in Figures \ref{fig:noisy} and \ref{fig:noisy_negative}, respectively. Since the exact (reference) solution is decreasing (cf. Figure \ref{fig:ref_sol2}), additive noise pushes the numerical solution upwards, causing it to decrease more slowly than the exact solution. Conversely, subtractive noise accelerates its decrease. In both cases, we consider the maximum possible error allowed by \eqref{eq:noise_function}.

Figures \ref{fig:noisy} and \ref{fig:noisy_negative} confirm the expected effect of the precision parameter: as the number of steps increases, the error contribution due to noise eventually starts to dominate the contribution due to the step size; cf. \eqref{eq:bound_th1}. Graphically, this manifests as a flattening of the error curve. For problem \eqref{eq:ivp3}, this effect becomes visible earlier for the randomized scheme and in the presence of subtractive noise.

Overall, numerical results illustrate the agreement between the observed and theoretical convergence rates established in this paper and in \cite{CzyzMorPrzyb}, and highlight the advantages of the randomized scheme.

\section{Conclusions and future work} \label{sec:conclusions}

In this paper, we have established $L^p(\Omega)$-error bounds for the randomized explicit Euler scheme applied to nonlinear ODEs under low regularity assumptions. In particular, we have considered a local one-sided Lipschitz condition and a local H\"older condition with respect to the space variable, as well as H\"older continuity in time. We have also allowed inexact evaluations of the right-hand side with additive noise satisfying an at-most-linear growth condition. The theoretical results capture the characteristic benefit of randomization, namely, a gain of $1/2$ in the convergence rate associated with temporal regularity.

Numerical experiments confirm the theoretical findings, with empirical convergence rates in line with the derived bounds. This supports the relevance of the proposed analysis for nonlinear problems under mild regularity assumptions.

Possible directions for future research include extending the analysis to higher-order methods, such as a randomized midpoint (RK2) method, under similarly mild assumptions. It would also be interesting to investigate lower bounds for the error and the optimality of the method in the information-based complexity setting. Finally, one could explore the construction of confidence intervals for the solution to quantify the uncertainty arising from randomized approximations and noisy evaluations.

{\bf Acknowledgements}
 This research was supported by the AGH University of Krakow, Poland, under grant no. 16.16.420.054, funded by the Polish Ministry of Science and Higher Education.

 {\bf Disclaimer}
 This paper was prepared for information purposes and is not a product of HSBC Bank Plc. or its affiliates. Neither HSBC Bank Plc. nor any of its affiliates make any explicit or implied representation or warranty and none of them accept any liability in connection with this paper, including, but not limited to, completeness, accuracy, reliability of information contained herein and the potential legal, compliance, tax or accounting effects thereof.

\end{document}